\documentclass{article}
\usepackage{amssymb}
\def\ba{\mathbf{a}}   
\def\bb{\mathbf{b}}   
\def\be{\mathbf{e}}   
\def\bbf{\mathbf{f}}   
\def\bff{\mathbf{f}}   
\def\bm{\mathbf{m}}   
\def\bn{\mathbf{n}}   
\def\bu{\mathbf{u}}   
\def\bv{\mathbf{v}}   
\def\bx{\mathbf{x}}   
\def\by{\mathbf{y}}   
\def\bz{\mathbf{z}}   

\def\bA{\mathbf{A}}   
\def\bB{\mathbf{B}}   

\def\C{\mathbb{C}} 
\def\G{\mathbb{G}} 

\def\R{\mathbb{R}}  

\def\no{\noindent}

\def\beq{\begin{equation}}
\def\eeq{\end{equation}}

\def\w{\wedge}

\def\no{\noindent}
\usepackage{helvet}
\usepackage{courier}
\usepackage{graphicx} 
\usepackage{array}    
\usepackage{tabularx}
\usepackage{amsfonts}
\usepackage{amssymb}
\usepackage{tikz-cd}
\def\bpm{\begin{pmatrix}}
\def\epm{\end{pmatrix}}
\begin{document}
	\title{What's in a Pauli Matrix?}
	\author{Garret Sobczyk
		\\ Universidad de las Am\'ericas-Puebla
		\\ Departamento de Actuar\'ia F\'isico-Matem\'aticas
		\\72820 Puebla, Pue., M\'exico
		\\ http://www.garretstar.com}
	\maketitle
	\begin{abstract} 
		 Why is it
		that after so many years matrices continue to play such an important roll in Physics and mathematics? Is there a geometric way of looking at matrices, and linear transformations in general, that lies at the roots of their success? We take an in depth look at the Pauli matrices, $2\times 2$ matrices over the complex numbers, and examine the various possible geometric interpretation of such matrices. The geometric interpretation of the Pauli matrices explored here natually extends to what the author has dubbed the study of {\it geometric matrices}. A geometric matrix is a matrix of order $2^n \times 2^n$ over the real or complex numbers, and has its geometric roots in its algebraically isomorphic Clifford geometric algebras, \cite{Hyprevisit}. 
			\smallskip
		\no {\em AMS Subject Classification:} 15A18, 15A66, 15B33, 83A05
		\smallskip
		
		\no {\em Keywords:} Clifford algebra, complex numbers, geometric algebra, idempotents, nilpotents, Pauli matrices.
		
	\end{abstract}
\newtheorem{thm3.1}{Lemma}
\newtheorem{thm3.2}{Theorem}
\newtheorem{thm3.3}[thm3.1]{Theorem}
\newtheorem{eigen}{Definition}
\newtheorem{specbasis}[thm3.2]{Theorem}
\newtheorem{thm3.4}[thm3.2]{Theorem}
\newtheorem{thm3.5}[thm3.2]{Theorem}

\section{The g-number system $\R(\ba, \bb  )$ }

	Sir Arthur Eddington (1882-1944) said, ``I cannot believe that anything so ugly as multiplication of matrices is an essential part of the scheme of nature", \cite{hays}.
	Tobias Danzig, in his {\it Number the Language of Science}, remarks ``These filing cabinets” are added and multiplied, and a whole calculus of
	matrices has been established which may be regarded as a continuation of the algebra of complex numbers", \cite{TD67}. The whole purpose of this article is to search for the magic of matrices in the most likely source of their power, their {\it geometric roots}. To make our task significantly easier, we restrict our attention to the simplest example of a powerful class of matrices, famously known as {\it Pauli matrices}. These matrices lie at the heart of the great discovery of quantum mechanics, which has revolutionized scientific and technological advancement over the last Century.

The real number system $\R $, which is at the heart of mathematics, is naturally pictured on the {\it real number line}. What makes the real number system so powerful is the ability to both {\it add} and {\it multiply} real numbers to get other real numbers. God saw that the real numbers were good, but the people were still not happy. Are there not more treasures to be found if we tinker just a bit with God's rules? For $r,s,t \in \R$, we have  
\begin{itemize}
	 \item[R1)] $rs =s r$ \quad {\it Commutative law of multiplication.}
	\item[R2)] $r(s+t) = rs + rt $ \quad {\it Distributive law of multiplication over addition.} 
	\item[R3)] $(rs)t = r(st)$ \quad {\it Associative law of multiplication.}
	\item[R4)]  $rs = 0 \ \ \iff r=0 \ \ {\rm and/or} \ \ s=0$.
\end{itemize}

 Let us introduce two {\it new numbers} $\ba$ and $\bb $ not in $\R $. To emphasize that
$\ba $ and $\bb$ are {\it not} real numbers, we write $\ba \notin \R$ and $\bb \notin \R $. To make the extended number system $\R(\ba , \bb )$, called {\it g-numbers} because of their forthcoming geometric interpretation, fully functional and compatible with $\R $, we extend the operations of addition and multiplication to include the new numbers $\ba, \bb$. This is accomplished by assuming that the extended numbers in $\R(\ba,\bb) $ obey exactly the same rules of addition and multiplication as do the numbers in $\R$, with the {\it exception} that $\ba \bb \ne \bb \ba $. Regarding the
the new numbers $\ba $ and $\bb $, they satisfy the following two special properties:
  \begin{itemize}
  	\item[N1)] $\ba^2 = 0 = \bb^2$. \quad {\it The new numbers $\ba $ and $\bb $ are called null-vectors or nilpotents}.
  	\item[N2)] $2 \ba \cdot \bb \equiv\ba \bb + \bb \ba =1$. \quad {\it The sum of $\ba \bb $ and $\bb \ba $ is $1$}.\footnote{By interpreting $\ba$ and $\bb$ as just being new kinds of vectors, it is natural to interpret $\ba \bb + \bb \ba$ as twice the inner product of these null-vectors.} 
    		
  \end{itemize} 
  
Clearly the non-zero g-numbers $\ba $ and $\bb $ cannot be real numbers because $\ba^2 = 0 = \bb^2$, since there are no non-zero real numbers with this property. Any g-number $g$ such that $g^2=0$ is said to be a {\it nilpotent}.
Also, the products $\ba \bb $ and $\bb \ba $ cannot be real numbers since
$\ba \bb \ne \bb \ba $. Never-the-less the property N2) tells us that the {\it sum}
$\ba \bb + \bb \ba =1\in \R$, providing a direct relationship between the extended g-numbers in $\R (\ba , \bb)$ and the real numbers, and showing that 
$\R \subset \R(\ba , \bb)$.

The {\it Multiplication Table} \ref{table1} for g-numbers is easily derived from the assumed properties N1) and N2), and the associative law. Since half of its entries are zeros, it is easily remembered. 
\begin{table}[h!]
	\begin{center}
		\caption{Multiplication table.}
		\label{table1}
		\begin{tabular}{c|c|c|c|c} 
			& $\ba $  & $ \bb $  & $\ba \bb $  & $ \bb \ba $ \cr 
			\hline
			$ \ba$  &   0   & $ \ba \bb $     & $  0 $         & $ \ba $ \cr 
			\hline         
			$ \bb $  &  $ \bb \ba $  & $ 0 $     & $  \bb $         & $ 0 $ \cr 
			\hline     
			$ \ba \bb $    &  $\ba $    &  0  &  $\ba \bb $    & 0  \cr
			\hline
			$ \bb \ba $  &  0   &  $\bb $   &  0     & $\bb \ba$   \cr
		\end{tabular}
	\end{center}
\end{table} 
To show that $\ba \bb \ba = \ba $, we use both properties N1) and N2), and in particular N2) to substitute in
$1 - \bb \ba $ for $\ba \bb$, getting
\[ \ba \bb \ba = (\ba \bb ) \ba = (1- \bb \ba)\ba = \ba - \bb \ba^2  = \ba , \]
and similarly, $\bb \ba \bb = \bb $. The same substitution works for showing that
\[ (\ba \bb )^2 = (1- \bb \ba  )(\ba \bb ) =\ba \bb +\bb \ba^2 \bb = \ba \bb ,  \]
and similarly that $(\bb \ba)^2= \bb \ba$. Any non-zero g-number $g$ with the property that $g^2=g$ is said to be an {\it idempotent}, so $\ba \bb $ and $\bb \ba $ are idempotents, and since 
\[ \ba \bb +\bb \ba =1, \quad {\rm and} \quad (\ba \bb)(\bb \ba)= 0 = (\bb \ba)  (\ba \bb), \]
they are said to {\it partition unity} and to be {\it mutually annihilating}, respectively. 

In addition, we assume that numbers in $\cal \R(\ba, \bb) $ {\it commute} with real numbers in $\R$, and that the associative and distributive properties R2) and R3) above remain valid for our new numbers. The g-numbers in the table are written as
a matrix,
\beq \R(\ba , \bb )=\pmatrix{\bb\ba  & \bb \cr \ba & \ba \bb }=\pmatrix{\bb \ba \cr  \ba} \pmatrix{ \bb  \ba & \bb}, \label{canbasis} \eeq 
and make up the {\it canonical nilpotent} or {\it null basis} of $\R(\ba , \bb ) $ over the real numbers.\footnote{Since $g$-numbers satisfy the same rules as do matrices,  matrices of $g$-numbers are well defined.} 
The last equality on the right expresses the nilpotent basis as the product of a column matrix of nilpotents with a row matrix of nilpotents. 
 
\section{Properties of g-numbers} 
 
Each g-number $g \in \R(\ba , \bb ) $ is uniquely specified by four real numbers.
In matrix form $[g]:=\pmatrix{g_{11} & g_{12} \cr g_{21} & g_{22} }$. Thus,
\beq  g = g_{11} \bb \ba + g_{12} \bb + g_{21}\ba + g_{22}\ba \bb =\pmatrix{\bb \ba &  \ba} [g] \pmatrix{ \bb  \ba \cr \bb} ,   \label{newnumbersN} \eeq
where $g_{11},g_{12},g_{21},g_{22}\in \R$.

  We can now easily derive the general rule for the addition and multiplication of two
 g-numbers $f, g \in \R(\ba,\bb) $. In addition to $g$, already defined, let
  \[  f = f_{11} \bb \ba + f_{12} \bb + f_{21}\ba + f_{22}\ba \bb=\pmatrix{\bb \ba &  \ba} [f] \pmatrix{ \bb  \ba \cr \bb}  \]
    Calculating $f+g$ and $f g$, we find that
\[ f+g = (f_{11}+g_{11})\bb \ba +(f_{12}+g_{12}) \bb +(f_{21}+g_{21})\ba +(f_{22}+g_{22})\ba \bb\] 
\beq =\pmatrix{\bb \ba &  \ba} \pmatrix{f_{11}+g_{11} & f_{12}+g_{12} \cr f_{21}+g_{21} & f_{22}+g_{22} } \pmatrix{ \bb  \ba \cr \bb}  =\pmatrix{\bb \ba &  \ba} ([f]+[g]) \pmatrix{ \bb  \ba \cr \bb} \label{matrixadd} \eeq
for addition, and
 \[ fg = (f_{11}g_{11}+f_{12}g_{21})\bb \ba + (f_{11}g_{12}+f_{12}g_{22}) \bb \]
  \[\quad \quad \quad \quad \quad \quad +(f_{21}g_{11}+f_{22}g_{21})\ba + (f_{21}g_{12}+f_{22}g_{22})\ba \bb\] 
 \[ =\pmatrix{\bb \ba &  \ba} \pmatrix{f_{11}g_{11}+f_{12}g_{21} & f_{11}g_{12}+f_{12}g_{22} \cr f_{21}g_{11}+f_{22}g_{21} &f_{21}g_{12}+f_{22}g_{22}}\pmatrix{ \bb  \ba \cr \bb}\] 
\beq =\pmatrix{\bb \ba &  \ba} ([f][g]) \pmatrix{ \bb  \ba \cr \bb} , \label{matrixmult} \eeq
for multiplication. Surprisingly, (\ref{matrixadd}) and (\ref{matrixmult}) reduce adding and multiplying $g$-numbers to the familiar rules for adding and multiplying $2 \times 2$ matrices,
  \[ [f+g]=[f]+[g]  \quad {\rm and} \quad [fg]=[f][g].  \]

To complete our new number system $\R(\ba,\bb)$, we define three powerful {\it conjugation operators} with respect to the canonical null basis (\ref{canbasis}). First note that each g-number $g\in \R(\ba,\bb) $ is the sum of 
two parts, $g = g_o+ g_e$, an {\it odd part $g_o$}, and an {\it even part} $g_e$, where 
\beq g_o:=  g_{12} \bb  + g_{21}\ba   , \ \ {\rm and} \ \ g_e := g_{11} \bb \ba  + g_{22}\ba \bb , \label{defoddeven}  \eeq
respectively. The odd part $g_o$ is a linear combination of the null vectors $\ba $ and $\bb $, and the even part is a linear combination of their products the idempotents $\bb \ba $ and $\ba \bb $.

 We define the {\it reverse} $g^\dagger$ of $g$, with respect to the null canonical basis (\ref{canbasis}), by
\beq g^\dagger := (g_o+g_e)^\dagger = g_o^\dagger + g_e^\dagger  = g_o + g_e^\dagger, \label{reverseop} \eeq
where
\[ g_o^\dagger := g_o \ \ {\rm and}  \ \ g_e^\dagger = 
             g_{11} \bb \ba  + g_{22}\ba \bb  .      \]
The reverse operation reverses the order of the multiplication of $\ba $ and $\bb $, {\it i.e.}, $(\ba \bb )^\dagger = \bb \ba $, leaving the odd part $g_o$ unaffected. It follows
that for $f,g \in \R(\ba,\bb) $,
\[  (f+g)^\dagger = f^\dagger + g^\dagger \ \ {\rm and} \ \ (fg)^\dagger = g^\dagger f^\dagger. \]

The {\it inversion} $g^-$ of $g$, with respect to the null canonical basis (\ref{canbasis}), is defined by 
\beq g^- := (g_o + g_e)^- = g_o^- + g_e^- =- g_o + g_e, \label{inversionop} \eeq
where
\[ g_o^- :=- g_{12}\bb   - g_{21}\ba = -g_o \ \ {\rm and}  \ \ g_e^- = g_e .      \]
The operation of inversion changes the {\it sign} of both $\ba $ and $\bb $, {\it i.e.}, $\ba^-=-\ba $ and $\bb^- = -\bb$, leaving $g_e$ unaffected. 
Clearly 
 \[  g_o = \frac{1}{2}(g-g^-), \quad {\rm and} \quad  
	g_e= \frac{1}{2}(g+g^-), \]
 and for $f,g \in \R(\ba,\bb)$,
\[  (f+g)^- = f^- + g^- \ \ {\rm and} \ \ (fg)^- = f^- g^-. \]

Combining the operations of reverse and inversion gives the third {\it mixed
conjugation}. For $g\in \R(\ba,\bb) $, the mixed conjugation $g^*$ of $g$, with respect to the standard canonical basis (\ref{canbasis}),  is defined by 
\beq g^* := (g^\dagger)^- = (g_o+g_e^\dagger)^-=-g_o +g_e^\dagger. \label{mixedop} \eeq 
The mixed conjugation of the sum and product of $f, g \in \R(\ba,\bb)$, satisfies 
\[ (f+g)^* = f^* + g^*= -(f_o +g_o)+(f_e^\dagger + g_e^\dagger), \]
and 
\[ (fg)^* = g^* f^* = (- g_o + g_e^\dagger)(- f_o + f_e^\dagger) =
(fg)^*_o+(fg)^*_e,    \]
where
\[(fg)^*_o = -(g_of_e^\dagger + g_e^\dagger f_o)\] 
and
   \[    (fg)^*_e =  g_o f_o +g_e^\dagger f_e^\dagger .  \]
   
 For $g\in \R(\ba,\bb)$, using the mixed conjugation, we define
\beq tr(g):= g+ g^* = g_e + g_e^\dagger = g_{11}+ g_{22}=tr[g]\in \R , \label{deftrace} \eeq
called the {\it trace} of $g$. 
 Also, using that an even g-number times an odd g-number is odd, with respect to the standard canonical basis (\ref{canbasis}), we calculate
\[ \det g := g g^* = (g_o+g_e) (-g_o+g_e^\dagger) \] 
\[=g_og_e^\dagger - (g_o g_e^\dagger)^\dagger + g_e g_e^\dagger - g_og_o \]  
\beq = g_e g_e^\dagger - g_og_o = g_{11}g_{22}-g_{12}g_{21}=\det[g] \in \R , \label{defdet} \eeq 
since
\[  g_e g_e^\dagger = (g_{11}\bb \ba + g_{22}\ba \bb ) (g_{11}\ba \bb + g_{22}\bb \ba ) = g_{11}g_{22} \]
and 
\[   g_o g_o = (g_{12}\bb  + g_{21}\ba)^2 = g_{12}g_{21}(\bb \ba + \ba \bb  ) = g_{12}g_{21}.\]

Given a g-number $g \in \R(\ba,\bb) $, when is there a $f\in \R(\ba,\bb)$ such that
$ g f = f g =1$? When such an $f$ exists, we say that $g^{-1}:= f$ is the {\it multiplicative inverse} of $g$. Since
\[ \frac{g^* g}{g^* g} = \frac{g g^*}{g g^*} =1 ,  \]
it immediately follows that 
\beq g^{-1} := \frac{g^*}{g g^*} =  \frac{g^*}{g^* g} = \frac{-g_o+g_e^\dagger}{g g^*} \label{inversegN} \eeq
provided that $g g^* \ne 0$. Whenever a g-number has the property that $\det g \ne 0$, $g$ is {\it non-singular}, and if $g g^*=0$, $g$ is {\it singular}. 
 
  Given g-numbers $f,g \in \R(\ba,\bb)$, the product $fg$ can be decomposed into even and odd parts with respect to the canonical null basis (\ref{canbasis}). We have 
  \[  fg = (fg)_o + (fg)_e\] 
  where
  \[ (fg)_o = f_o g_e + f_e g_o \ \ {\rm and}  \ \ (fg)_e = f_og_o +  f_e g_e  \]   
  
The product of two g-numbers $fg$ can also be decomposed into the sum of a {\it symmetric part} $f\circ g$ and a {\it skew-symmetric part} $f\otimes g$. We have
\beq fg = \frac{1}{2}(fg + gf) +  \frac{1}{2}(fg - gf) = f\circ g + f \otimes g,
           \label{symskewsym} \eeq
  where $f\circ g := \frac{1}{2}(fg + gf)$ and $f \otimes g :=\frac{1}{2}(fg - gf)$.
    For the null vectors $\ba $ and $\bb $, we find that
   \beq \ba \circ \bb  \equiv \ba \cdot \bb= \frac{1}{2}(\ba \bb + \bb \ba)=\frac{1}{2} \ \ {\rm and} \ \
   \ba \otimes \bb \equiv\ba \w \bb= \frac{1}{2}(\ba \bb- \bb \ba). \label{relinout} \eeq 
   Squaring $\ba \otimes \bb = \ba \w \bb$, gives
   \[ (\ba \w \bb )^2= \frac{1}{4}(\ba \bb- \bb \ba)^2= \frac{1}{4}\Big((\ba \bb)^2+ (\bb \ba)^2\Big) = \frac{1}{4},  \]
   so $  (\ba \w \bb )^2 =( \ba \cdot \bb)^2=\frac{1}{4} $. 
  It follows that the idempotents
  \[  \ba \bb = \ba \cdot \bb + \ba \w \bb = \frac{1}{2} (1+2 \ba \w \bb) \ \ {\rm and} \ \   \bb \ba = \bb \cdot \ba + \bb \w \ba = \frac{1}{2} (1-2 \ba \w \bb).  \]

  \section{Geometry of $\R(\ba,\bb)$} 
  
 \begin{figure}[h]
 	\centering
 	\includegraphics[width=0.45\linewidth]{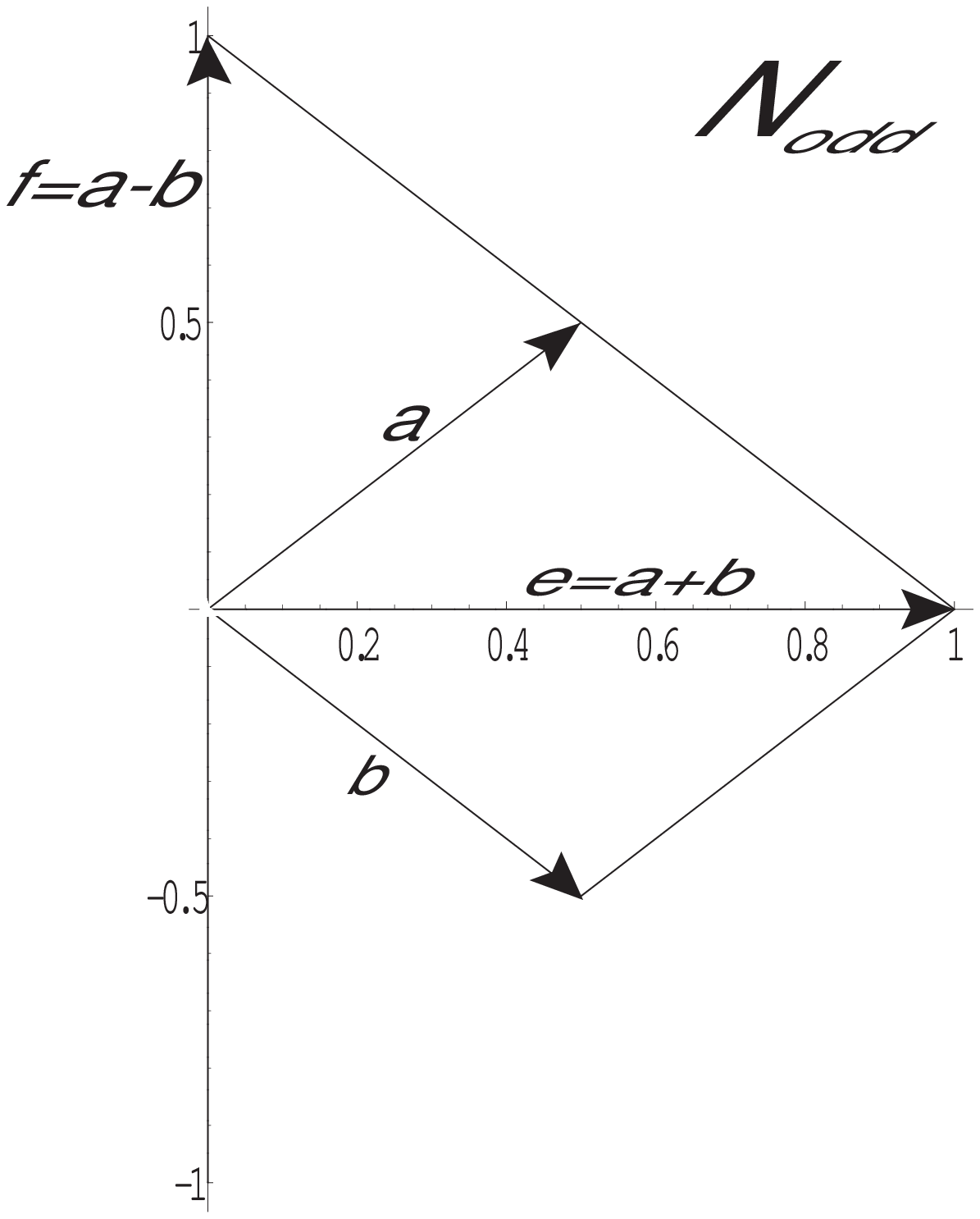}
 	\includegraphics[width=0.45\linewidth]{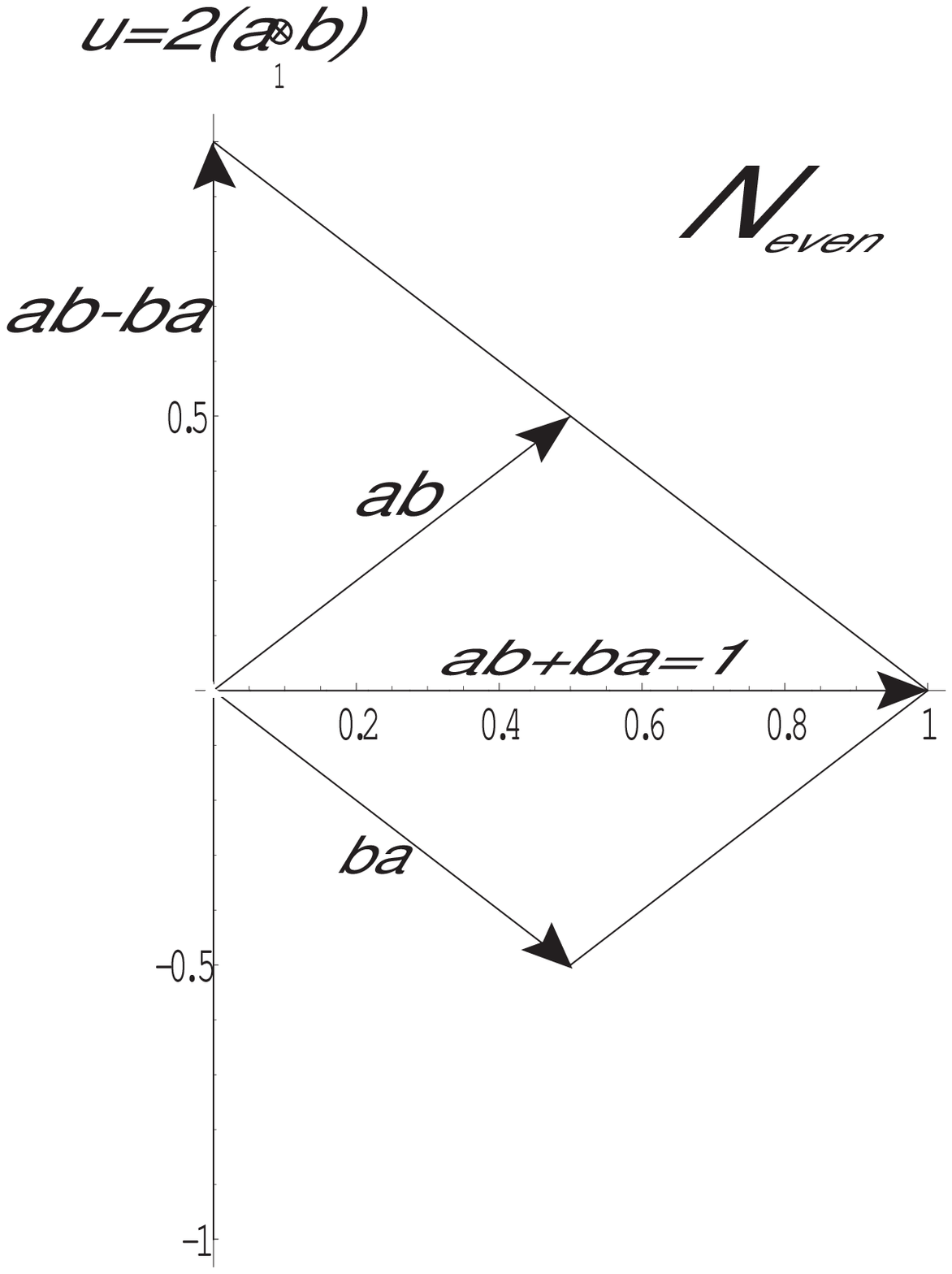}
 	\caption{The odd number plane $\R(\ba,\bb) _o$.\quad The even number plane $\R(\ba,\bb) _e$.}
 	\label{evenoddhyp}
 \end{figure}

  Much conceptual clarity is gained when it is possible to pictorially represent
  fundamental concepts. We picture the odd and even parts of a g-number
  \[ g=g_{11}\bb \ba +g_{12}\bb + g_{21}\ba + g_{22}\ba \bb = g_o+g_e \in \R(\ba,\bb) \]
  separately in the odd and even g-number planes $\R(\ba,\bb)_o$ and $\R(\ba,\bb)_e$, respectively,
Figure \ref{evenoddhyp}. Because $\ba $ and $\bb $ are nilpotents, they are pictured on a 2-dimensional {\it null-cone} in $\R(\ba,\bb)_o$. 
 Similarly, since $\ba \bb$ and $\bb \ba $ are {\it singular idempotents}, they define the $2$-dimensional null-cone in $\R(\ba,\bb)_e$. 

Notices the different conventions that have been adopted in Figure \ref{evenoddhyp}. In the odd $g$-number plane, the vector $\be:=\ba+\bb$ lies along the positive $x$-axis, and
 vector $\bff := \ba - \bb$ lies along the positive $y$-axis. On the other hand, in the even $g$-number plane $1 = \ba \bb + \bb \ba$ and $u := \ba \bb- \bb \ba $ are chosen along $x$- and $y$-axes, respectively, to agree with the conventions established for the equivalent hyperbolic number plane discussed in \cite{S1}. 

The orthonormal unit vectors $\be=\ba + \bb $ and $\bff=\ba -\bb $ make up the {\it standard basis} of $\R^{1,1}$, and generate the real geometric algebra 
\beq  \G_{1,1}:=\G(\R^{1,1}) =\R(\be, \bff)\equiv \R(\ba, \bb).  \label{stanbasis11} \eeq  
 It is easy verify that the unit vectors
$\be $ and $\bff$ satisfy the basic rules
\beq \be^2= 1 = - \bff^2, \ {\rm and} \ \be \bff = -\bff \be. \label{stdbasis11} \eeq   
Notice that the unit bivector 
\beq \bbf \be =(\ba - \bb)(\ba + \bb)=2 \ba \w \bb ,\label{bivectoru} \eeq 
unlike a unit Euclidean bivector, has square plus one 
\[  ( \bbf \be)^2 = \bbf \be \bff \be = -\bff^2 \be^2 =1.   \] 

In terms of the standard basis, the g-number
(\ref{newnumbersN}) has the form
 \[ g =
 \frac{1}{2}\Big[tr(g)+(g_{12}+g_{21})(\ba +\bb)+(g_{22}-g_{11})2\ba \w \bb +(g_{21}-g_{12})(\ba - \bb )
 \Big]  \]
\beq = \alpha_0 + \alpha_1 \be  + \alpha_2 \bff \be  +\alpha_3 \bbf = \alpha_0 + \bv , \label{ginstdb} \eeq 
for
\[ \alpha_0:=\frac{1}{2}(g_{22}+g_{11}), \ \alpha_1:=\frac{1}{2}(g_{21}+g_{12}), \] \[\alpha_2:=\frac{1}{2}(g_{22}-g_{11}), \ \alpha_3:=\frac{1}{2}(g_{21}-g_{12}).       \]
A general g-number $g\in \G_{1,1}$ consists of two parts, a {\it scalar part}
$\alpha_0:=\frac{1}{2}tr(g) $, and a {\it vector part} $\bv:= \alpha_1 \be  + \alpha_2 \bff \be  +\alpha_3 \bbf $. The vector part $\bv$ of $g$ is a misnomer because it is a linear combination of not only $\be $ and $\bff$, but also of the bivector $\bff \be $. As will be become increasingly clear, the concept of a vector and a bivector are frame-related, or observer-dependent quantities. What is important here is that the basis elements $\be, \bff \be , \bff$ are mutually anti-commutative. 

In the standard basis, the determinant of $g$ takes the form
\[ \det g= gg^* =\alpha_0^2 -\alpha_1^2-\alpha_2^2+\alpha_3^2 =\alpha_0^2-\bv^2.   \]
 A g-number $g$ is said to be {\it hyperbolic}, {\it parabolic} or {\it Euclidean} if  
\beq \big(g-\alpha_0)^2= \bv^2  \begin{array}{cr}>0,\\=0,\\<0,
\end{array}  \label{hypparaeuclid} \eeq
respectively. If $g$ is non-singular and hyperbolic, then $g$ has one of the {\it hyperbolic  Euler forms}  
\beq g= \pm \rho e^{\phi u} \quad {\rm or} \quad g= \pm \rho u e^{\phi u},  \label{hypeuler} \eeq
for $u=\frac{\alpha_1 \be  + \alpha_2 2\ba \w \bb +\alpha_3 \bbf}{|\alpha_1 \be  + \alpha_2 2\ba \w \bb +\alpha_3 \bbf|}$ and $\pm \rho$ chosen appropriately, \cite{S1}. If $g$ is parabolic and $gg^* =0$, then $g$ is a nilpotent. If $g$ is parabolic and $gg^* \ne 0$, then $g$ has the Euler form
\beq g=\alpha_0 e^{\frac{\alpha}{\alpha_0} (\alpha_1 \be  + \alpha_2 2\ba \w \bb +\alpha_3 \bbf)}=\alpha_0 +\alpha \bn ,\label{lighteuler} \eeq
for $\bn = \alpha_1 \be  + \alpha_2 \bff \be +\alpha_3 \bbf $.
 In the case that $g=\alpha_0 + \bv $ is Euclidean,
 then $g$ has the Euclidean  Euler form
 \beq g = re^{i \theta},    \label{euclideuler} \eeq
 for $i=\frac{\alpha_1 \be  + \alpha_2 2\ba \w \bb +\alpha_3 \bbf}{|\alpha_1 \be  + \alpha_2 2\ba \w \bb +\alpha_3 \bbf|}$ and $r:=\sqrt{g g*}=\sqrt{\alpha_0^2+\bv^2}$.
  Once the Euler form of a g-number is found, and $\det g\ne 0$, it is easy to find the
 Euler form of $g^{-1}$. For example, if
 $g= \rho u e^{\phi u}$, then
 $g^{-1}= \frac{1}{\rho} u e^{-\phi u}$.  

The real number system $\R $ has no non-zero nilpotents. We have defined the extended number system $\R(\ba,\bb)$ in terms of new nilpotents $\ba$ and $\bb$ which satisfy the multiplication Table \ref{table1}, and the rules N1) and N2). The question arises: What is the most general number in $\bn\in\R(\ba,\bb)$ which is a non-zero nilpotent? Let $[\bn]:=\pmatrix{n_1 & n_2 \cr n_3 & n_4}$ be the real number matrix of $n$. Then
\beq [\bn^2]=\pmatrix{n_1 & n_2 \cr n_3 & n_4}^2=\pmatrix{n_1^2+n_2 n_3 & n_1n_2+n_2 n_4 \cr n_1n_3+n_3 n_4 & n_2 n_3+n_4^2}=0, \label{realnilpotent}  \eeq 
or equivalently,
\[ \bn^2 = (n_1^2+n_2 n_3)\bb \ba + ( n_1n_2+n_2 n_4)\bb +( n_1n_3+n_3 n_4)\ba + (n_2 n_3+n_4^2)\ba \bb =0.  \]

The real solutions to (\ref{realnilpotent}) can be broken into two cases:

\begin{itemize}
	\item[Case 1.]  $n_2 n_3 = 0 \ \ \iff \ \ n_1=n_4=0$ giving non-trivial solutions:
	$n_2\ne 0 \ \ n_3=0$ and $n_3\ne 0 \ \ n_2 =0$. 
\item[Case 2.] $n_2 n_3 \ne 0 \ \ \iff \ \ n_2\ne 0, \ n_3 \ne 0\ {\rm and} \ n_1 = - n_4$. \end{itemize} 
It follows that all non-trivial real solutions of (\ref{realnilpotent}) are of the form
\beq \bn = n_2 \bb + n_3 \ba + 2 n_4 \ba \w \bb   \label{nilpotentcan} \eeq
where $n_2 n_3 + n_4^2 = 0$. When $n_4=0$, we fall back to Case 1 with only two non-trivial solutions. If $n_4 \ne 0$, then any non-trivial solution must satisfy the
condition that $n_2 n_3 < 0$. In terms of the standard canonical basis (\ref{stanbasis11}), (\ref{stdbasis11}),  any non-trivial solution
can be put in the form
\beq \bn = \alpha_1 \be+ \alpha_2 \bff \be  + \alpha_3 \bff  ,    \label{nilpotentstdcan} \eeq 
where $\alpha_1^2+\alpha_2^2 - \alpha_3^2=0$ for $\alpha_1,\alpha_2,\alpha_3\in \R $. Recalling 
$\bbf \be=2 \ba \w \bb $, we find that
$\alpha_1=\frac{n_3+n_2}{2}$, $\alpha_2=\frac{n_3-n_2}{2}$ and $\alpha_3 = n_4 = -n_1$ relating (\ref{nilpotentcan}) and (\ref{nilpotentstdcan}). The set of all null g-numbers $\bn$ of the form (\ref{realnilpotent}) or (\ref{nilpotentcan}) defines the {\it null cone} ${\cal N} \in \G_{1,1} $, shown in Figure \ref{hypparaeuclid}. 

\begin{figure}[h]
	\centering
	\includegraphics[width=0.45\linewidth]{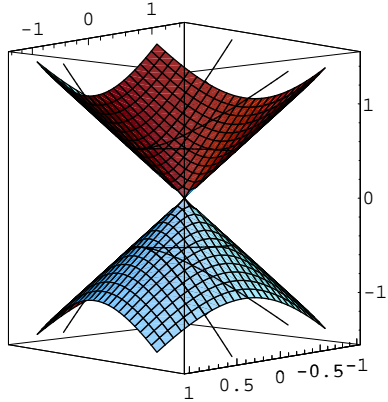}
	\caption{Hyperbolic, Parabolic and Euclidean mappings of the null cone $\cal N$ in $\G_{1,1}$. }
	\label{hypparaeuclid}
	\end{figure}

Given the nilpotent $\ba$ and a general nilpotent $\bn$ such that $[\bn] =\pmatrix{-n_4 & n_2 \cr n_3 & n_4}$, where $n_4^2 +n_2 n_3=0$, the following theorem shows that there is always a g-number $g\in \G_{1,1}$ such that $\bn = g \ba g^{-1}$.

\begin{thm3.4} Given a nilpotent $\bn$ with the matrix $[\bn ]$ and $n_2\ne 0$, the g-number $g=\ba + n_2 \bb + n_4 \ba \bb$, with the matrix $[g]:=\pmatrix{0 & n_2 \cr 1 & n_4}$, has the property that $\bn =g \ba g^{-1}$. \label{thm3.4}
	\end{thm3.4}
	
\no {\bf Proof.} 	
	 
 We wish to find a non-singular $g$-number $g$ such that $g\ba g^{-1}=\bn $, or equivalently, 
\[ [g][\ba][g^{*}]=\det[g] [\bn]  .  \]
Calculating,
\[ \pmatrix{g_{11}& g_{12}\cr g_{21} & g_{22}} \pmatrix{0 & 0 \cr 1 & 0}\pmatrix{g_{22}& -g_{12}\cr -g_{21} & g_{11}}=\pmatrix{g_{12}& 0\cr g_{22} & 0} \pmatrix{g_{22} & -g_{12} \cr -g_{21} & g_{11}}    \]
\[=\pmatrix{g_{12}g_{22}& -g_{12}^2\cr g_{22}^2 & -g_{12}g_{22}}= \det[g]\pmatrix{-n_{4}& n_{2}\cr n_{3} & n_{4}},  \] 
or 
\[n_4=-\frac{g_{12}g_{22}}{\det[g]} , \quad n_2=-\frac{g_{12}^2}{\det[g]} , \quad  n_3=\frac{g_{22}^2}{\det[g]}, \]
implying that $n_2$ and $n_3$ must have different signs. The choice $g_{21}=1, g_{22}=n_4$ and $g_{12}=n_2$, giving $\det g = -n_2$, solves these conditions and completes the proof.

\hfill $\square$

{\bf Corollary} Given a nilpotent $\bn$ with the matrix $[\bn ]$ and $n_2\ne 0$, the g-number $g=\frac{n_4}{2} + \bv$, with $\bv = \frac{1+n_2}{2}\be +\frac{n_4}{2}\bff \be + \frac{1-n_2}{2} \bff$ has the property that $\bn = g \ba g^{-1}$.
\smallskip

 The Corollary is merely a restatement of the Theorem in terms of the standard basis
(\ref{stdbasis11}). Thus, the inner automorphism defined by $g$ is hyperbolic, parabolic, or Euclidean (\ref{hypparaeuclid}), according to whether 
\[ \bv^2 = \frac{1}{4}(n_4^2+ 4 n_2) \begin{array}{cr}>0,\\=0,\\<0,
\end{array}  \]
 respectively.
 
 As an example, consider nilpotent $\bn= 4 \ba \w \bb + \bb - 4 \ba$ with the matrix $[\bn]=\pmatrix{-2& 1 \cr -4 & 2}$. 
The g-number $g= 1 + \ba + \bb + 2 \ba \w \bb$ with the matrix $[g]=\pmatrix{0&1 \cr 1 & 2}$ is hyperbolic and has the property $g \ba g^{-1} = \bn$,
since $ ( \ba + \bb + 2 \ba \w \bb)^2 =\frac{3}{2}$. 
The same hyperbolic mapping applied to $\bb$ gives
$ \ba = g\bb g^{-1}$, or equivalently, $g^{-1} \ba g = \bb$ .

As a further example, let $ h=2 \ba \w \bb + \bb - 2 \ba$ so that $[h]= \pmatrix{-1 & 1 \cr -2 & 1}$, and let 
\[  [\bn] = [h\ba h^{-1}]=\pmatrix{1 & -1 \cr 1 & -1}, \  {\rm so \ that} \ \bn = -2 \ba \w \bb + \ba -\bb,   \]
and 
\[  [\bm] = [h\bb h^{-1}]=\pmatrix{-2 & 1 \cr -4 &2}, \ {\rm so \ that} \  \bm =4 \ba \w \bb + \bb - 4 \ba.    \]
Since $(2 \ba \w \bb +\bb - 2 \ba)^2 = -1$, $h$ defines a Euclidean rotation. It is interesting to compare this example to the previous example.

The following Lemma summarizes the results derived in (\ref{hypeuler}), (\ref{lighteuler}), and (\ref{euclideuler}).
\begin{thm3.1}Let $g\in \R(\ba ,\bb) $ be a non-singular geometric number. Then $g$ has one of the three Euler forms, according to whether $g$ is hyperbolic, parabolic, or Euclidean. 
	
	a) If $g$ is hyperbolic, then 
	\[ g=\pm \rho e^{\phi u} \quad {\rm or} \quad g=\pm \rho ue^{\phi u},  \]
	for appropriately chosen $u$, $\phi$, and $\rho$. 
	
	b) If $g$ is parabolic, then 
	\[ g = \alpha_0 e^{\frac{\alpha}{\alpha_0}\bn}=\alpha_0+\alpha \bn ,   \]
	with $\alpha$ and $\bn$ chosen as in (\ref{lighteuler}).
	
	c) If $g$ is Euclidean, then
	\[ g= r e^{i \theta} ,  \]
where $\theta$, $i$, and $r$ are chosen as in (\ref{euclideuler}).	

\end{thm3.1}
 
 Let $\G_{1,1}=\R(\ba , \bb)=\R(\be,\bbf )$ be the geometric algebra defined by the canonical null vector basis $\{\ba ,\bb\} $. The following theorem characterizes the {\it regrading} of the geometric algebra induced by non-singular g-number $g \in \G_{1,1}$ taking the null vector basis $\{\ba ,\bb\} $ into a {\it relative null vector basis}
 $\{\bA ,\bB\} $, where $\bA := g\ba g^{-1}$ and $\bB := g\bb g^{-1}$. In Figure \ref{hypparaeuclid}, hyperbolic, parabolic, and Euclidean mappings of $\bA$ and $\bB$, as a function of $\theta$, $\alpha$ and $\phi$, defined by the respective Euler forms, are shown.    
  
  \begin{thm3.2} Each non-singular hyperbolic, parabolic, or Euclidean, g-number $g \in \G_{1,1}$ induces a regrading of $\G_{1,1}$ into $\G_{1,1}^\prime := \R(\bA , \bB)$,
  	with the relative cononical null vector basis    
 	\beq {\cal B} :=\pmatrix{\bB \bA & \bB \cr \bA & \bA \bB }  \label{canbasisAB} \eeq
satisfying 
\[ \bA^2=0, \  \bB^2=0, \   \bA \circ \bB =\frac{1}{2}, \] 
where
 \[ \pmatrix{\bB \bA & \bB \cr \bA & \bA \bB }:=g\pmatrix{\bb \ba & \bb \cr \ba & \ba \bb}g^{-1} .  \]
   \end{thm3.2}
     
  {\bf Proof:} A non-singular g-number is either hyperbolic, parabolic, or Euclidean, having the Euler forms given in (\ref{hypeuler}), (\ref{lighteuler}), or (\ref{euclideuler}), respectively. In each of these cases it is clear that
   \[ \pmatrix{\bB \bA & \bB \cr \bA & \bA \bB }:=\pmatrix{g\bb g^{-1}g\ba g & g\bb g^{-1}\cr g\ba g^{-1} & g\ba g^{-1}g \bb g^{-1}}=g\pmatrix{\bb \ba & \bb \cr \ba & \ba \bb}g^{-1} .  \]

Theorem 2 shows that partitioning a geometric number $g=g_e + g_o$ into {\it even} (scalar and bivector)
and {\it odd} (vector) parts, satisfying the multiplication rules given in Table \ref{table1}, and with an relative inner and outer product (\ref{relinout}), is a {\it relative} concept with the new null vectors
$\bA $ and $\bB $ defined by 
   \[ \bA = g\ba g^{-1}, \ {\rm and} \ \bB = g\bb g^{-1}.    \]
 Whereas the algebras $\R(\ba,\bb)\widetilde{=}\R(\bA, \bB)$ are isomorphic, the nilpotents (vectors) in $\R(\bA,\bB)$ are a mixture of vectors and bivectors in $\R(\ba , \bb)$. Never-the-less, any such partition defines the (relative) geometric algebra $\G_{1,1}$ as given in  (\ref{stanbasis11}). Each non-singular g-number defines a partitioning of the elements of $\G_{1,1}$ into
 relative null vectors on the null cone $\cal N$.
 \hfill $\square$

 The most general idempotent $\bB \bA \in\G_{1,1}$ has the form
 \beq \bB \bA = \frac{1}{2}\Big( 1+ \alpha_1  \be+ \alpha_2  \bff \be+ \alpha_3  \bff  \Big),
   \label{genidempotent11} \eeq  
  where $\alpha_3=\pm \sqrt{\alpha_1^2 + \alpha_2^2 -1}$ and $\alpha_1,\alpha_2,\alpha_3 \in \R $. Given the idempotent
  $\bB \bA$ the g-number defined by the matrix
  \[ [g]:= \pmatrix{ 1-\alpha_2 & \alpha_1 - \alpha_3 \cr 1+\alpha_2 & -\alpha_1+\alpha_3},   \]
   has the property that $g\bB \bA g^{-1}= \bb \ba $ provided that $\det g = 2(\alpha_3-\alpha_1)\ne 0$. In this case, the matrices of $\bA := g^{-1}\ba g $ and $\bB:=g^{-1}\bb g  $ are given by
  \beq   [\bA] :=\pmatrix{\frac{1-\alpha_2}{2} & \frac{\alpha_1-\alpha_3}{2} \cr \frac{(1 -\alpha_2)^2}{2(\alpha_3-\alpha_1)} & - \frac{1-\alpha_2}{2}}, \ \  [\bB] :=\pmatrix{\frac{1+\alpha_2}{2} & \frac{\alpha_3-\alpha_1}{2} \cr \frac{-(1+\alpha_2)^2}{2(\alpha_3-\alpha_1)} & - \frac{1+\alpha_2}{2}} \label{nilpotentsAB11}. \eeq

 The geometric algebras $\G_{1,1}$ and $\G_{2,0}$ are algebraically isomorphic. The isomorphism $f: \G_{1,1}\widetilde = \G_{2,0}$, specified by the mapping
 \beq \be \leftrightarrow \be_1, \ \bbf \leftrightarrow \be_{12}, \ \be \bbf  \leftrightarrow \be_2,      \label{g11tog2}  \eeq
  is not given by an inner-automorphism defined by a non-singular $g \in \G_{1,1}$. We have 
 \[ \G_{1,1}:=\R( \be,\bff)=span_\R \{1,\be_1,\bff_1,\be_1\bff_1 \} \widetilde = \] 
 \beq  span_\R \{1,\be_1,\be_2,\be_{12} \}= \R(\be_1,\be_2)=:\G_{2} . \label{geos2} \eeq

   \section{Structure of a geometric number}
   
   The fact that the canonical null basis (\ref{canbasis}) consists only of g-numbers which are nilpotents, or a product of nilpotents, suggests that these g-numbers are of fundamental importance. Traditionally, in linear algebra the characteristic polynomial plays a crucial role. The equivalent property of a g-number is given in the definition below.   
  
  \begin{eigen} Given a geometric number $g=\alpha_0 + \bv \in \R(\ba,\bb)=\G_{1,1}  $. The characteristic
  	polynomial of $g$ is 
  	\beq \varphi_g(\lambda):=(\lambda-\alpha_0)^2-\bv^2. \label{charpoly} \eeq
  	 The real or complex roots $\lambda=\alpha_0 \pm v$ of this	polynomial, for $v = \sqrt{\bv^2}$,  are the eigenvalues of $g$.  
        	\end{eigen}
      The structure of the geometric number $g$ is completely determined by its characteristic polynomial $\varphi_g(\lambda)$, \cite{S2,S0,S3}. The eigenvalues of $g$ can be
      either real or complex numbers. If $g$ is hyperbolic or parabolic, then the eigenvalues of $g$ are real. On the other hand, if $g$ is Euclidean, then the eigenvalues are conjugate complex. Special attention is given to the case of complex eigenvalues. Complex eigenvalues are formally assumed to commute with g-numbers in $\G_{1,1}$. 
      
 The different canonical forms of a g-number in $ \G_{1,1}$ are given in
  \begin{thm3.4} A g-number $g=\alpha_0 + \bv \in \G_{1,1} $ has one of the three canonical forms:  
  	\begin{itemize}
 		\item[i)] a) If $g$ is hyperbolic, then for $\rho = \sqrt{\bv^2 }>0$,
 \beq g =\lambda_1 \hat \bv_+ + \lambda_2 \hat \bv_-                 \label{hypgnum} \eeq	

 where $\lambda_1:= \alpha_0+ \rho$, $\lambda_2:= \alpha_0- \rho $ and
 \[ \hat \bv_+ := \frac{1}{2}(1+ \hat \bv) \quad {\rm and} \quad \hat \bv_- := \frac{1}{2}(1- \hat \bv).   \]
 \item[(i)] b) When $\bv = 0$, then $g=\alpha_0$ and $\lambda_1 = \lambda_2 = \alpha_0$.
 		\item[ii)] If $g$ is Euclidean, then for $i:=\sqrt{-1}$ and $\rho =\sqrt{-\bv^2}$, 
 		\beq g =\lambda_1 \hat \bv_+ + \lambda_2 \hat \bv_-                 \label{euclidgnum} \eeq		
 		where $\lambda_1:= \alpha_0- \rho i $, $\lambda_2:= \alpha_0+ \rho i$, and
 		\[ \hat \bv_+ := \frac{1}{2}(1+ i\hat \bv) \quad {\rm and} \quad \hat \bv_- := \frac{1}{2}(1-i \hat \bv).   \]
 		\item[iii)] If $g$ is parabolic and $\bv \ne 0$, then
\beq  g= \alpha_0 + \bn, \label{paragnum} \eeq
where $\bn = \bv$ is a nilpotent.

 	\end{itemize}
 	  \end{thm3.4}
       
     {\bf Proof:} The proof, a straight forward verification, is omitted. 
        
        \hfill $\square$
        
           Given that $f\in \G_{1,1}$ is type {\it i)}, so that $f=\lambda_1 \bv_+ +\lambda_2 \bv_-$, by multiplying both sides of this equation on the right by $\hat \bv_+$ and $\hat \bv_-$, successively, we get
           \beq  f \bv_+ = \lambda_1 \bv_+, \ \ {\rm and} \ \ f \bv_- = \lambda_2 \bv_-, \label{eigenpotentf} \eeq 
           respectively. We say that $\bv_+$ and $\bv_-$ are {\it eigenpotents} for the
           respective eigenvalues $\lambda_1$ and $\lambda_2$. When $f=\lambda + \bn$ for type {\it ii)}, multiplying on the right by $\bn$ gives 
           $  f\bn = \lambda \bn$. In this case, we also say that $\bn$ is an {\it eigenpotent} of $f$. It is interesting that for type {\it i}) $f\in \G_{1,1}$, that the eigenpotents are idempotents, whereas for type {\it ii}) $f$, the eigenpotent is a nilpotent. More fundamentally, the following theorem shows each g-number defines a relative rest-frame of eigenpotents which are always nilpotents.

 \begin{specbasis} i) If $f= \alpha_0 + \bv$ is type i), so that
 	 \[ f=\lambda_1 \bv_+ +\lambda_2 \bv_-,  \] 
 then there exits nilpotents $\bA, \bB \in \cal N$, such that $\bA \circ \bB =\frac{1}{2}$, which are eigenpotents of $f$ satisfying
 \[ f \bB = \lambda_1 \bB \quad {\rm and} \quad f\bA = \lambda_1 \bA .     \]
 
 \no ii) If $f$ is type ii), so that $f = \alpha_0 + \bn $, then there exists a relative
 canonical null basis $\cal B$ such that the matrix of $f$ has the form
 \[ [f]_{\cal B}=\pmatrix{\alpha_0 & 1 \cr 0 & \alpha_0}.    \]

 \end{specbasis}	 	
  
  {\bf Proof:} i) Applying (\ref{genidempotent11}) and (\ref{nilpotentsAB11}) to the nonzero singular idempotent $\bB \bA :=\hat \bv_+$, we can find a non-singular $g\in \R(\ba,\bb) $, and construct a relative canonical null basis (\ref{canbasis}), such that
   \[ {\cal B} := \pmatrix{\bB \bA & \bB \cr \bA & \bA \bB }=g^{-1} \pmatrix{\bb \ba & \bb \cr \ba & \ba \bb }g,   \]
   where $\hat \bv_+=\bB \bA$, $\hat \bv_- = \bA \bB $, $\bA \circ \bB = \frac{1}{2}$, and $\bA^2=0 = \bB^2$, so that
      \[ f=\lambda_1 \bB \bA + \lambda_2 \bA \bB .      \]
  Multiplying both sides this equation on the right by $\bB $, and then by $\bA $, gives
  \beq  f\bB = \lambda_1 \bB \bA \bB = \lambda_1 \bB \ \ {\rm and}  \ \
      f\bA = \lambda_2 \bA \bB \bA = \lambda_2 \bA, \label{eigennilps}  \eeq
  respectively. Equations (\ref{eigenpotentf}) and (\ref{eigennilps}) are
  equivalent since we can easily get back the first equation from the second. In the canonical null basis $\cal B$, the matrix of $f$ is 
  \[ [f]_{\cal B} = \pmatrix{\lambda_1 & 0 \cr 0 & \lambda_2}. \] 
  
           {\it ii}) When $f= \lambda + \bn$ for the nilpotent $\bn$, we use {\bf Theorem \ref{thm3.4}} 
  to find a g-number $g\in \G_{1,1} $ such that    
    \beq \bn= g^{-1}\bb g .\label{compatiblegm}  \eeq 
 The relative canonical basis $\cal B$ of $g$ is then defined by
     \[  {\cal B} := \pmatrix{\bB \bA  & \bB  \cr \bA  & \bA \bB}  =
     g^{-1} \pmatrix{\bb \ba & \bb \cr \ba & \ba \bb }g,  \]
  where 
  \[ \bB:=g^{-1} \bb  g = \bn \quad {\rm and} \quad \bA:=g^{-1}\ba g = \bm .  \]
  With respect to this relative basis ${\cal B}$, the matrix of $f$ is 
  \[ [f]_{\cal B}:= \pmatrix{\lambda & 1 \cr 0 & \lambda}.  \] 
  
  \hfill $\square$

  There are a number of vector analysis like identities that are useful when carrying out calculations with the vector parts of g-numbers. Let 
   \[\bx := (x_1,x_2, x_3 )\pmatrix{\be \cr \bbf \be \cr \bff }, \ 
   \by := (y_1,y_2, y_3 )\pmatrix{\be \cr \bbf \be \cr \bff }, \ 
   \bz := (z_1,z_2,z_3 )\pmatrix{\be \cr \bbf \be \cr \bff }. \]
   The symmetric or {\it scalar product} of $\bx $ and $\by $ is
   \beq \bx \circ \by = x_1 y_1+x_2 y_2 - x_3 y_3.  \label{scalarprodxy} \eeq
   The anti-symmetric or {\it cross product} of $\bx$ and $\by$ is
 \beq \bx \otimes \by =\det \pmatrix{\be & \bff \be & -\bff \cr
 	                      x_1 & x_2 & x_3 \cr
 	                      y_1 & y_2 & y_3 }         .  \label{crossprodxy} \eeq
 In addition, there are two triple product,
   \beq \bx \circ (\by \otimes \bz) =  \det \pmatrix{x_1 & x_2  & x_3 \cr
   	y_1 & y_2 & y_3 \cr
   	z_1 & z_2 & z_3 }    =( \bx \otimes \by) \circ \bz  ,  \label{tribleprodxyz1} \eeq	                    
  and
   \beq \bx \otimes ( \by \otimes \bz) =(\bx \circ \by) \bz - (\bx \circ \bz )\by      .  \label{tripleprodxyz2} \eeq 
  The proofs of these formulas is left to the reader.

  \section{Geometric algebras of $2\times 2$ matrices}
 In previous sections, we have seen how real $2\times 2$ matrices are the {\it coordinates} of g-numbers in $\R(\ba ,\bb )$, or in the corresponding geometric algebra $\G_{1,1}$. Geometric algebras assign geometric meaning to what otherwise are just a tables of numbers,
 \cite{S08,Sob2012}. However, in studying the structure of real g-numbers, the embarrassment of complex eigenvalues arises. Just as the real number system $\R$ is extended to the complex number system $\C $, real g-numbers $\R(\ba,\bb) $ are extended to the complex g-numbers $\C(\ba,\bb)$. The real and complex g-numbers $\R(\ba,\bb) $ and $\C(\ba,\bb)$ are algebraically isomorphic to the Clifford geometric algebras $\G_{1,1}$ and $\G_{1,2}$, respectively. 
  
  In terms of its matrix $[f]$, since $(\bb \ba)(\bb \ba)=\bb \ba $ and $\bb \ba [f]=[f]\bb \ba $, 
  \[  f   =\pmatrix{\bb \ba & \ba}[f]\pmatrix{\bb \ba \cr \bb }  = \pmatrix{\ba \bb & \bb } [f^\dagger]^T\pmatrix{\bb \ba  \cr \bb}\]
    \beq
  =  \pmatrix{\bb \ba & \ba} \pmatrix{f_{11} & f_{12} \cr f_{21} & f_{22}}\pmatrix{\bb \ba \cr \bb } =f_{11} \bb \ba + f_{12}\bb + f_{21}\ba + f_{22} \ba \bb.  \label{matrixoff}  \eeq
  Furthermore,
  \[ f^\dagger = \pmatrix{\ba \bb  & \bb}[f]^T\pmatrix{\ba \bb \cr \ba  }, \ \ 
  {\rm and} \ \  f^* = \pmatrix{\ba \bb  & -\bb}[f]^T\pmatrix{\ba \bb \cr -\ba  }, \]
  where $[f]^T$ is the {\it transpose} of the matrix $[f]$.
   
  The equation (\ref{matrixoff}) can be directly solved for the matrix $[f]$ of $f$.
  Multiplying equation (\ref{matrixoff}) on the left and right by
  $\pmatrix{\bb \ba \cr \bb}$ and $\pmatrix{ \bb \ba  & \ba}$, respectively, gives the  equation
    \[ \pmatrix{\bb \ba \cr \bb}f\pmatrix{ \bb \ba  & \ba}   = \pmatrix{\bb \ba \cr \bb}\pmatrix{\bb \ba & \ba}[f]\pmatrix{\bb \ba \cr \bb } \pmatrix{ \bb \ba  & \ba}\]
    \[ = \pmatrix{\bb \ba & 0 \cr 0 & \bb  \ba}[f]\pmatrix{\bb \ba & 0 \cr 0 & \bb \ba }
   = \bb \ba [f].   \]
  Similarly, multiplying equation (\ref{matrixoff}) on the left and right by
  $\pmatrix{  \ba  \cr \ba \bb }$ and $\pmatrix{\bb  & \ba \bb}$, respectively, gives 
  \[ \pmatrix{ \ba \cr \ba \bb}f\pmatrix{ \bb   & \ba \bb}   = \pmatrix{ \ba \cr \ba \bb}\pmatrix{\bb \ba & \ba}[f]\pmatrix{\bb \ba  \cr  \bb } \pmatrix{ \bb   & \ba \bb}\]
  \[ = \pmatrix{ \ba & 0 \cr 0 &   \ba}[f]\pmatrix{\bb  & 0 \cr 0 & \bb  }
  = \ba \bb [f].   \]
 Adding these two equations together give the desired result
 \[ [f]= \pmatrix{\bb \ba \cr \bb}f\pmatrix{ \bb \ba  & \ba} + \pmatrix{ \ba \cr \ba \bb}f\pmatrix{ \bb   & \ba \bb} \] 
 \beq =  \pmatrix{\bb \ba f \bb \ba + \ba f \bb & 
 	 \bb \ba f \ba + \ba f \ba \bb \cr
 	 \bb f \bb \ba + \ba \bb f \bb &\bb f \ba + \ba \bb f \ba \bb }  .  \label{matrixf} \eeq

  The geometric algebra $  \G_{1,1}$ is defined by
  \[  \G_{1,1}:=\R(\be,\bbf)   \]
  where $\be^2 = 1 = - \bbf^2 $, and $\be \bbf = -\bbf \be .$ 
  The geometric algebra $\G_{1,1}$ is the real number system $\R $ extended to include the new anticommuting square roots $\be , \bbf$ of $\pm 1$, respectively. Since by (\ref{stdbasis11})
  \[ \pmatrix{\be & \bbf}=\pmatrix{ \ba & \bb}\pmatrix{1 & 1 \cr 1 & -1} \quad \iff  \quad  \pmatrix{ \ba & \bb} = \frac{1}{2}\pmatrix{\be & \bbf} \pmatrix{1 & 1 \cr 1 & -1} ,  \]
  it follows that
  \[ \G_{1,1} := \R(\be , \bbf )= \R(\ba , \bb),     \] 
 reflects only a change of basis of $\G_{1,1}$. Recall also the case of the real geometric
 algebra $\G_{2,0}$, obtained in (\ref{geos2}) by re-interpreting the elements of $\G_{1,1}$,
 \[  \be \to \be_1, \ \ \be \bff \to \be_2, \ \ {\rm and} \ \ \bbf \to \be_{12}.     \]
 Whereas the algebras $\G_{1,1}$ and $\G_{2,0}$ are algebraically isomorphic, there is no inner automorphism relating them, since such an inner automorphism would violate the famous {\it Law of inertia}, \cite[p.297,334]{Gant1960}.

 We have defined $\R(\ba ,\bb)$ to be the real number system extended to include the new elements $\ba, \bb$. Because of the problem of {\it complex} eigenvalues, we further extend $\R(\be,\bff) $ to $\C(\be,\bff)$ by allowing
the matrix $[g]_\C$ to consists of {\it complex numbers}. Thus 
\[  g = \pmatrix{\bb \ba & \ba } [g]_\C \pmatrix{\bb \ba \cr \bb} , \]
where $[g]_c $ is $2\times 2$ matrix over the complex numbers $\C $. Of course, we have make the additional adhoc assumption that complex numbers commute with $\ba $ and $\bb $, and therefore with all the g-numbers in $\R(\ba,\bb) $.

On the level of geometric algebras, we can give complex g-numbers in $\C(\be,\bff  )$  different
 geometric interpretations. The geometric algebras $\G_{1,2}$ and $\G_3$ are defined
 by 
 \beq \G_{1,2} := \R(\be_1, \bbf_1, \bbf_2), \ {\rm and} \ \G_3:=\R(\be_1,\be_2,\be_3)     \label{geoalg12g13}  \eeq
where $\{\be_1, \be_2, \be_3\}$, and $\{  \bbf_1, \bbf_2\}$, are new {\it anti-commuting}
square roots of $\pm 1$, respectively. 

We have
 \beq \G_{1,2}:= \R(\be_1, \bbf_1, \bbf_2)\widetilde{=} \C(\be, \bff  )  \label{defg12} \eeq 
 for
$\be_1:= \be$, $\bff_1:= \bff$, and $\bff_2:= i \be_1 \bff_1$, and
\beq \G_3:=\R(\be_1,\be_2,\be_3) \widetilde = \C(\be, \bff ),  \label{defg30} \eeq
for $\be_1:= \be, \be_2:= i \bff$ and $ \be_3:=\be \bff $. In both of these
real geometric algebras $\G_{1,2}$ and $\G_3$, the formally adhoc imaginary number $i:=\sqrt{-1}$ takes on the geometric interpretation of the unit pseudoscalar or volume element in these respective algebras. This follows directly from the calculations
\[  \be_1 \bff_1 \bff_2 = \be \bff (i \be \bff) = i, \ {\rm and} \ 
    \be_1 e_2 \be_3  = \be (i\bff ) (\be \bff) = i. \] 
  In $\G_3$, using (\ref{defg30}) the famous Pauli matrices of $\be_1,\be_2,\be_3$ are easily found to be
  \beq [\be_1]=\pmatrix{0& 1 \cr 1 & 0} , [\be_2]=\pmatrix{0& -i \cr i & 0}, [\be_3 ]=\pmatrix{1& 0 \cr 0 & -1}. \label{paulimat} \eeq

A complex g-number $f$ in $\G_{1,2}$, or $\G_3$, is defined by its complex matrix $[f]_\C$ in the canonical null basis (\ref{newnumbersN}), by 
\[ f=f_{11}\bb \ba + f_{12} \bb + f_{21} \ba + f_{22} \ba \bb = \alpha_0 + \bv,\]
where $f_{jk}\in \C\,   \widetilde =\,  \G_{1,2}^{0+3}$, or in $\G_{3}^{0+3}$, respectively. In the standard complex basis, the $\alpha_0$ and $\bv=\alpha_1 \be + \alpha_2 \bff \be + \alpha_3 \bff$ are defined in the same way as in the standard real basis (\ref{ginstdb}). Each of
the complex scalars $f_{jk} $ are defined by
\[  f_{jk}= x_{jk} + iy_{jk}, \quad {\rm where} \quad i:=\be_1 \bbf_1 \bbf_2 \ {\rm or} \ i=\be_{123},       \]
respectively, for $x_{jk},y_{jk} \in \R$, with complex conjugates 
\[    \overline{f_{ij}}:=f_{ij}^-= f_{ij}^\dagger .\] 
Geometric algebras exist for arbitrary signatures,
\[ \G_{p,q}:=\G(\R^{p,q} )= \R(\be_1, \ldots, \be_p, \bff_1,\ldots, \bff_q), \]
\cite{geoS2017}.   

As an example, the most general Hermitian g-number w.r.t the basis of Pauli vectors
(\ref{paulimat}) is $h=\alpha_0+\bu$ for 
\[ \bu:=u_1 \be_1 +u_2 \be_2+ u_3 \be_3 = \rho \hat \bu   \]
where $u_1,u_2,u_3 \in \R$ and  $\rho:= \sqrt{u_1^2+u_2^2+u_3^2}$. The eigenvalues and eigenpotents of $h$ are 
$\lambda_\pm := \alpha_0 \pm \rho$, and 
$\frac{1}{2}(1\pm \hat \bu )$, respectively. Choosing $\bB \bA =\frac{1}{2}(1+ \hat \bu )$,
the g-number whose matrix is
\[  [g] = \pmatrix{\frac{u_1+ i u_2}{\rho - u_3} & 1 \cr -\frac{u_1+ i u_2}{\rho + u_3} &1  } ,  \]
has the property that $\bB \bA = g^{-1} \bb \ba g$. It follows that
\[ h= \lambda_+ \bB \bA + \lambda_- \bA \bB .  \]
Using that $\bA = g^{-1}\ba g$, and $\bB = g^{-1}\bb g$, we find the eigenpotents 
\[  \bA := \frac{\rho + u_3}{2 \rho}\pmatrix{-1 & - \frac{u_1-i u_2}{\rho + u_3} \cr \frac{u_1+i u_2}{\rho - u_3} & 1} \ \ {\rm and} \ \  \bB := \frac{\rho - u_3}{2 \rho}\pmatrix{-1 &  \frac{u_1-i u_2}{\rho - u_3} \cr -\frac{u_1+i u_2}{\rho + u_3} & 1},   \]
which satisfy $h \bA = \lambda_- \bA$ and $h \bB := \lambda_+ \bB$. A {\it spacetime} vector analysis was developed in \cite{S81}. See \cite{Sob2012,SNF,{S1/2},Shopf2015} for many other applications.

\end{document}